\documentclass[leqno]{article}

\begin{document}

\title{Potpourri, 6}

\author{Stephen William Semmes	\\
	Rice University		\\
	Houston, Texas}

\date{}

\maketitle

	Let $E$ be a nonempty finite set, and let $V$ be a real or
complex vector space.  In these notes we write $\mathcal{F}(E, V)$ for
the vector space of $V$-valued functions on $E$.

	Suppose that $E_1$, $E_2$ are nonempty finite sets, and that
$V$ is the vector space of real or complex-valued functions on $E_2$.
In this case we can identify $\mathcal{F}(E_1, V)$ with the vector
space of real or complex-valued functions on the Cartesian product
$E_1 \times E_2$, as appropriate.

	In any event $\mathcal{F}(E, V)$ is basically the same as the
tensor product of the vector space of real or complex-valued functions
on $E$, as appropriate, with $V$.  Of course the tensor product of two
vector spaces, both real or both complex, is symmetric up to canonical
isomorphism in the two vector spaces.  With $\mathcal{F}(E, V)$ we may
unravel this symmetry a bit, depending on the circumstances.

	For each $x \in E$, let $\delta_x(y)$ denote the scalar-valued
function on $E$ which is equal to $1$ when $y = x$ and to $0$ when $y
\ne x$.  These functions form a basis for the scalar-valued functions
on $E$, and in effect we may treat the vector space of scalar-valued
functions on $E$ as having a distinguished basis in this way.

	We can also think of the scalar-valued functions on $E$ as
being a commutative algebra over the real or complex numbers, as
appropriate, using ordinary pointwise multiplication as the product in
the algebra.  We can think of $\mathcal{F}(E, V)$ as being a module
over the commutative algebra of scalar-valued functions, since we can
multiply a vector-valued function by a scalar-valued function.

	Suppose that $A$ is any linear transformation on $V$.  There
is a unique linear transformation $\widetilde{A}$ on $\mathcal{F}(E,
V)$ such that $\widetilde{A}(f)(x) = A(f(x))$ for every $V$-valued
function $f$ on $E$ and every $x \in E$.  Alternatively,
$\widetilde{A}$ is characterized by the property that $\widetilde{A}(f
\, v) = f \, A(v)$ for every scalar-valued function $f$ on $E$ and
every vector $v \in V$.

	Now suppose that $T$ is any linear transformation acting on
$\mathcal{F}(E, {\bf R})$ or $\mathcal{F}(E, {\bf C})$, as
appropriate.  There is a unique linear transformation $\widehat{T}$ on
$\mathcal{F}(E, V)$ such that $\widehat{T}(f \, v) = T(f) \, v$ for
all scalar-valued functions $f$ on $E$ and all $v \in V$.

	Thus one can also think of $\mathcal{F}(E, V)$ as modules over
the algebras of linear transformations acting on the vector space of
scalar-valued functions on $E$ and on $V$.  The actions of
$\widetilde{A}$ and $\widehat{T}$ on $\mathcal{F}(E, V)$ described in
the preceding paragraphs automatically commute, by construction.

	Let us focus for the moment on scalar-valued functions.  Let
$E$ be a nonempty finite set, and let $w(x)$ be a weight on $E$.  To
be more precise, $w(x)$ is a real-valued function on $E$ such that
$w(x) > 0$ for all $x \in E$.

	If $f(x)$ is a real or complex-valued function on $E$,
and $p$ is a positive real number, put
\begin{equation}
	\|f\|_{w, p} = \bigg(\sum_{x \in E} |f(x)|^p \, w(x) \bigg)^{1/p}.
\end{equation}
When $w(x) = 1$ for all $x \in E$ we may simply write $\|f\|_p$ for
this quantity.

	For $p = \infty$ put
\begin{equation}
	\|f\|_\infty = \|f\|_{w, \infty} = \max \{|f(x)| : x \in E\}.
\end{equation}
Thus the weight plays no role when $p = \infty$.

	Suppose that $w \equiv 1$.  Clearly
\begin{equation}
	\|f\|_\infty \le \|f\|_p
\end{equation}
for all real or complex-valued functions $f$ on $E$ and all positive
real numbers $p$.  Using this one can check that
\begin{equation}
	\|f\|_q \le \|f\|_p
\end{equation}
when $0 < p \le q \le \infty$.

	Now suppose that $w$ defines a probability distribution on $E$,
which is to say that
\begin{equation}
	\sum_{x \in E} w(x) = 1.
\end{equation}
In this event
\begin{equation}
	\|f\|_{w, p} \le \|f\|_\infty
\end{equation}
for all real or complex-valued functions $f$ on $E$ and all positive
real numbers $p$.  If $r$ is a real number, $r \ge 1$, and $\phi(x)$
is a nonnegative real-valued function on $E$, then
\begin{equation}
	\bigg(\sum_{x \in E} \phi(x) \, w(x) \bigg)^r 
		\le \sum_{x \in E} \phi(x)^r \, w(x).
\end{equation}
by the convexity of the function $t^r$ on the nonnegative real
numbers.  It follows that
\begin{equation}
	\|f\|_{w, p} \le \|f\|_{w, q}
\end{equation}
when $0 < p \le q \le \infty$.

	Let $V$ be a real or complex vector space, and let $N(v)$ be a
real-valued function on $V$ such that $N(v) \ge 0$ for all $v \in V$,
$N(v) = 0$ if and only if $v = 0$, and $N(\alpha \, v) = |\alpha| \,
N(v)$ for all scalars $\alpha$ and all vectors $v \in V$.  If
furthermore
\begin{equation}
	N(v + w) \le N(v) + N(w)
\end{equation}
for all $v, w \in V$, then we say that $N$ defines a norm on $V$.  In
the presence of the other conditions, this is equivalent to saying
that the closed unit ball associated to $N$,
\begin{equation}
	\{v \in V : N(v) \le 1\},
\end{equation}
is a convex subset of $V$.  If there is a positive real number $C$
such that
\begin{equation}
	N(v + w) \le C \, (N(v) + N(w))
\end{equation}
for all $v, w \in V$, then we say that $N$ defines a quasinorm on $V$.
A special case of this occurs when there is a positive real number
$p$, $p \le 1$, such that
\begin{equation}
	N(v + w)^p \le N(v)^p + N(w)^p
\end{equation}
for all $v, w \in V$.

	For all $p$, $0 < p \le \infty$, $\|f\|_{w, p}$ satisifies the
positivity and homogeneity conditions mentioned at the beginning of
the previous paragraph.  When $p = 1$ or $p = \infty$ one can easily
verify from the definitions that $\|f\|_{w, p}$ satisifies the
triangle inequality and therefore defines a norm on the vector space
of real or complex-valued functions on $E$.  The triangle inequality
also holds when $1 < p < \infty$, as one can check by showing that the
associated closed unit ball is convex using the convexity of $t^p$ on
the nonnegative real numbers.  When $0 < p \le 1$ one has that
\begin{equation}
	\|f_1 + f_2\|_{w, p}^p \le \|f_1\|_{w, p}^p + \|f_2\|_{w, p}^p
\end{equation}
for all real or complex-valued functions $f_1$, $f_2$ on $E$.  This
follows from the fact that
\begin{equation}
	(a + b)^p \le a^p + b^p
\end{equation}
for all nonnegative real numbers $a$, $b$.

	Suppose that $f_1$, $f_2$ are real or complex-valued functions
on $E$, $1 \le p, q \le \infty$, and 
\begin{equation}
	\frac{1}{p} + \frac{1}{q} = 1.
\end{equation}
If either $p$ or $q$ is infinite, then its reciprocal is equal to $0$.
H\"older's inequality states that
\begin{equation}
	\biggl|\sum_{x \in E} f_1(x) \, f_2(x) \, w(x) \biggr|
		\le \|f_1\|_{w, p} \, \|f_2\|_{w, q}
\end{equation}
in this case.  Of course one might as well restrict one's attention to
functions $f_1$, $f_2$ which are real-valued and nonnegative.  If $p
=1$ and $q = \infty$, or $q = 1$ and $p = \infty$, then H\"older's
inequality can be verified directly.

	Suppose that $1 < p, q < \infty$.  If $a$, $b$ are nonnegative
real numbers, then
\begin{equation}
	a \, b \le \frac{a^p}{p} + \frac{b^q}{q}.
\end{equation}
This can be seen as a special case of the convexity of the exponential
function, for instance.  Hence
\begin{equation}
	|f_1(x) \, f_2(x)| \le \frac{|f_1(x)|^p}{p} + \frac{|f_2(x)|^q}{q}
\end{equation}
for all $x \in E$, and therefore
\begin{equation}
	\biggl|\sum_{x \in E} f_1(x) \, f_2(x) \, w(x) \biggr|
		\le \frac{\|f_1\|_{w, p}^p}{p} + \frac{\|f_2\|_{w, q}^q}{q}.
\end{equation}
This implies H\"older's inequality when $\|f_1\|_{w, p}$ and
$\|f_2\|_{w, q}$ are equal to $1$, and the general case follows from a
scaling argument.

	If $f_1$, $f_2$ are functions on $E$ and $0 < p, q, r \le \infty$,
\begin{equation}
	\frac{1}{r} = \frac{1}{p} + \frac{1}{q},
\end{equation}
then
\begin{equation}
	\|f_1 \, f_2\|_{w, r} \le \|f_1\|_{w, p} \, \|f_2\|_{w, q},
\end{equation}
by H\"older's inequality, with the obvious conventions for infinite
exponents.  For a single function $f$ on $E$ we have that
\begin{equation}
	\|f\|_{w, r} \le \|f\|_{w, p}^{p/r} \, \|f\|_{w, q}^{q/r},
\end{equation}
again by H\"older's inequality.  One can also use H\"older's
inequality to check that $\|f\|_{w, p}$ satisfies the triangle
inequality when $1 < p < \infty$, also known as Minkowski's
inequality.  Namely,
\begin{eqnarray}
\lefteqn{\qquad \sum_{x \in E} |f_1(x) + f_2(x)|^p \, w(x)}	\\
	&&	\le \sum_{x \in E} |f_1(x)| \, |f_1(x) + f_2(x)|^{p-1} \, w(x)
	   + \sum_{x \in E} |f_2(x)| \, |f_1(x) + f_2(x)|^{p-1} \, w(x),
						\nonumber
\end{eqnarray}
and one can apply H\"older's inequality to each of the sums on the
right side.

	Let $h$ be a real or complex-valued function on $E$, and
define a linear transformation $M_h$ on the vector space of real or
complex-valued functions on $E$, as appropriate, by
\begin{equation}
	M_h(f)(x) = h(x) \, f(x),
\end{equation}
so that $M_h(f)$ is the product of $h$ and $f$.  If $p$, $q$, $r$ are
as in the preceding paragraph, then we get that
\begin{equation}
	\|M_h(f)\|_{w, r} \le \|h\|_{w, q} \, \|f\|_{w, p}
\end{equation}
for all functions $f$ on $E$.  This inequality is sharp in the sense
that for each nonzero function $h$ on $E$ there is a nonzero function
$f$ on $E$ such that equality holds and indeed such that 
\begin{equation}
	|h(x) \, f(x)|^r = |h(x)|^q = |f(x)|^p
\end{equation}
for all $x \in E$, at least when $p, q, r < \infty$.  For infinite
exponents one can make simple adjustments to this last part.

	Let $T$ be a linear transformation on the vector space of real
or complex-valued functions on $E$, and fix $r$, $s$ with $0 < r, s
\le \infty$.  Suppose that $k$ is a nonnegative real number such that
\begin{equation}
	\|T(f)\|_{w, s} \le k \, \|f\|_{w, r}
\end{equation}
for all functions $f$ on $E$.  Let $h_1$, $h_2$ be real or
complex-valued functions on $E$, as appropriate, and suppose that $0 <
q_1, q_2 \le \infty$ and that $q_1 \ge r$.  Choose $p$, $t$ with $0 <
p, t \le \infty$ such that
\begin{equation}
	\frac{1}{r} = \frac{1}{p} + \frac{1}{q_1}
\end{equation}
and
\begin{equation}
	\frac{1}{t} = \frac{1}{q_2} + \frac{1}{s}.
\end{equation}
It follows from the earlier inequalities that
\begin{equation}
	\|(M_{h_2} \circ T \circ M_{h_1})(\phi)\|_{w, t}
	   \le k \, \|h_1\|_{w, q_1} \, \|h_2\|_{w, q_2} \, \|\phi\|_{w, p}
\end{equation}
for all functions $\phi$ on $E$.  Conversely, suppose that $k$ is a
nonnegative real number such that this inequality holds for all
functions $h_1$, $h_2$, and $\phi$ on $E$.  Then one can check that
$\|T(f)\|_{w, s} \le k \, \|f\|_{w, r}$ for all functions $f$ on $E$.
This is basically the same as the sharpness of the inequality for the
norm of $M_h(f)$ discussed in the previous paragraph.

	Now suppose that $V$ is a real or complex vector space with
quasinorm norm $\|\cdot \|_V$.  If $0 < p \le \infty$ and $f$ is a
$V$-valued function on $E$, put
\begin{equation}
   \|f\|_{w, p, V} = \bigg(\sum_{x \in E} \|f(x)\|_V^p \, w(x) \bigg)^{1/p}
\end{equation}
when $p < \infty$ and
\begin{equation}
	\|f\|_{\infty, V} = \|f\|_{w, \infty, V}
		= \max \{\|f(x)\|_V : x \in E\}.
\end{equation}
This is the same as $\|\cdot \|_{w, p}$ applied to the real-valued
function $\|f(x)\|_V$ on $E$.

	One can check that $\|f\|_{w, p, V}$ is a quasinorm on
$\mathcal{F}(E, V)$ for all $p$.  If $\|v\|_V$ is a norm on $V$ and $1
\le p \le \infty$, then $\|f\|_{w, p, V}$ is a norm on $\mathcal{F}(E,
V)$.  If $0 < p \le 1$ and $\|v\|_V$ satisfies
\begin{equation}
	\|v + w\|_V^p \le \|v\|_V^p + \|w\|_V^p
\end{equation}
for all $v, w \in V$, then $\|f\|_{w, p, V}$ satisfies the analogous
condition on $\mathcal{F}(E, V)$.  Observe that if $\|v\|_V$ satisfies
this condition and $0 < q \le p$, then we also have that
\begin{equation}
	\|v + w\|_V^q \le \|v\|_V^q + \|w\|_V^q
\end{equation}
for all $v, w \in V$.

	Let $h$ be a real or complex-valued function on $E$, as
appropriate, and let us again write $M_h$ for the linear operator on
$\mathcal{F}(E, V)$ given by multiplication by $h$.  Thus $M_h(f)(x) =
h(x) \, f(x)$ for all $V$-valued functions $f$ on $E$.
If $0 < p, q, r \le \infty$ and $1/r = 1/p + 1/q$, then
\begin{equation}
	\|M_h(f)\|_{w, r, V} \le \|h\|_{w, q} \, \|f\|_{w, p, V}
\end{equation}
for all $V$-valued functions $f$ on $E$.

	As in the scalar case, suppose that $0 < q_1, q_2, r, s \le
\infty$ and $q_1 \ge r$, and choose $p$, $t$ with $0 < p, t \le
\infty$ so that $1/r = 1/p + 1/q_1$ and $1/t = 1/q_2 + 1/s$.  If $Z$
is a linear transformation on $\mathcal{F}(E, V)$ and $k$ is a
nonnegative real number, then
\begin{equation}
	\|Z(f)\|_{w, s, V} \le k \, \|f\|_{w, r, V}
\end{equation}
for all $V$-valued functions $f$ on $E$ if and only if
\begin{equation}
	\|(M_{h_2} \circ Z \circ M_{h_1})(\phi)\|_{w, t, Z}
	   \le k \, \|h_1\|_{w, q_1} \, \|h_2\|_{w, q_2} \, \|\phi\|_{w, s, V}
\end{equation}
for all real or complex-valued functions $h_1$, $h_2$ on $E$, as
appropriate, and all $V$-valued functions $\phi$ on $E$.  In
particular this applies to the case where $Z = \widehat{T}$ for some
linear transformation $T$ acting on the vector space of real or
complex-valued functions on $E$, as appropriate.  In this event
$M_{h_2} \circ Z \circ M_{h_1}$, as a linear transformation on
$V$-valued functions, is the same as 
\begin{equation}
	\widehat{R}, \ R = M_{h_2} \circ T \circ M_{h_1},
\end{equation}
where $M_{h_1}$, $M_{h_2}$ refer to the multiplication operators on
scalar-valued functions in the definition of $R$.

	Let $a(x, y)$ be a real or complex-valued function on $E
\times E$.  This leads to a linear transformation $A$ acting on real
or complex-valued functions on $E$, as appropriate, given by
\begin{equation}
	A(f)(x) = \sum_{y \in E} a(x, y) \, f(y) \, w(y),
\end{equation}
and every such linear transformation can be expressed in this manner.
Consider the condition
\begin{equation}
	\|A(f)\|_{w, 1} \le \|f\|_{w, \infty}
\end{equation}
for all functions $f$ on $E$.  This holds in particular if
\begin{equation}
	\sum_{x \in E} \sum_{y \in E} |a(x, y)| \, w(x) \, w(y) \le 1.
\end{equation}

	Suppose that $T$ is a linear transformation acting on real or
complex-valued functions on $E$.  Let $r$, $s$ be exponents with $1
\le r, s \le \infty$, and assume that
\begin{equation}
	\|T(f)\|_{w, s} \le \|f\|_{w, r}
\end{equation}
for all functions $f$ on $E$.  Let $q$, $1 \le q \le \infty$, be the
exponent conjugate to $s$, so that $1/q + 1/s = 1$, and let $h_1$,
$h_2$ be real or complex-valued functions on $E$, as appropriate,
such that
\begin{equation}
	\|h_1\|_{w, r}, \|h_2\|_{w, q} \le 1.
\end{equation}
Under these conditions, the linear transformation
\begin{equation}
	A = M_{h_2} \circ T \circ M_{h_1}
\end{equation}
acting on functions on $E$ satisfies the conditions described in the
previous paragraph.

	Let $W$ be a finite-dimensional real or complex vector space
equipped with a norm $\|\cdot \|_W$.  Consider expressions for linear
transformations on $W$ of the form
\begin{equation}
	\sum_{\ell=1}^n c_\ell \, \lambda_\ell \, w_\ell,
\end{equation}
where the coefficients $c_\ell$ are real or complex numbers, as
appropriate, the $w_\ell$'s are elements of $W$ with $\|w_\ell\|_W \le
1$ for all $\ell$, and the $\lambda_\ell$'s are linear functionals on
$W$ such that $|\lambda_\ell(u)| \le \|u\|_W$ for all $u \in W$, which
is to say that the $\lambda_\ell$'s have dual norm less than or equal
to $1$.  Every linear transformation on $W$ can be expressed as a sum
of this kind.

	By definition, the trace norm of a linear transformation on
$W$ is equal to the infimum of
\begin{equation}
	\sum_{\ell = 1}^n |\alpha_\ell|
\end{equation}
over all representations of the linear transformation as in the
previous paragraph.  Let $p$ be a positive real number such that $p
\le 1$.  We can consider more generally the infimum of
\begin{equation}
	\bigg( \sum_{\ell = 1}^n |\alpha_\ell|^p \bigg)^{1/p}
\end{equation}
over all such representations of a given linear transformation.  This
quantity is greater than or equal to the trace norm, and is
monotonically increasing as $p$ decreases.  It defines a quasinorm on
the vector space of linear transformations on $W$ which satisfies the
generalization of the triangle inequality in which the $p$th power of
the quasinorm is subadditive.

\end{document}